\documentclass[12pt]{article}
 \usepackage{color}

\usepackage{latexsym}
\usepackage{epic,eepic,epsfig}
\usepackage{amssymb}
\usepackage{amsfonts,srcltx}

\newtheorem{theorem}{Theorem}[section]
\newtheorem{prop}[theorem]{Proposition}
\newtheorem{lem}[theorem]{Lemma}
\newtheorem{rem}[theorem]{Remark}

\newtheorem{defi}[theorem]{Definition}
\newtheorem{exam}[theorem]{Example}

\newcommand{\mod}{\, \hbox{mod} \, }

\newcommand{\Aut}{\hbox{Aut}\, }

\newcommand{\'}{\' \i}

\newcommand{\qed}{\hfill \mbox{\raisebox{0.7ex}{\fbox{}}}}

\def\demo{{\bf Proof}\hskip10pt}

 \def\MM{{\cal M}}

\def\o1{\overline 1} 

\def\o{\overline}     

\def\di{\bigm|} \def\lg{\langle} \def\rg{\rangle}
\def\nd{\mathrel{\bigm|\kern-.7em/}}

\def\f{\noindent}

\def\PSL{\hbox{\rm PSL}}

\def\Aut{\hbox{\rm Aut}}

\def\Syl{\hbox{\rm Syl}}

\def\Mon{\hbox{\rm Mon}}
\def\Aut{\hbox{\rm Aut}}
\def\mod{\hbox{\rm mod }}

\def\char{\hbox{\rm \,char\,}}
\newcommand{\B}[1]{\mathbb #1}

\begin{document}

\baselineskip=16pt
\newcounter{teocounter}
\addtocounter{teocounter}{1}
\newcounter{defcounter}
\addtocounter{defcounter}{1}
\newcounter{figcounter}
\addtocounter{figcounter}{1}
\newcounter{itc}
\addtocounter{itc}{1}
\newcounter{Lorrms}
\setcounter{Lorrms}{2} \vspace{2.0 cm}

\begin{center}

{\Large   Orientably-Regular $\pi$-Maps  and  Regular $\pi$-Maps}
\end{center}

\begin{center}  Xiaogang Li and Yao Tian$^{*}$
\end{center}

\begin{center}
{\footnotesize  School of Mathematical Sciences\\ Capital Normal University \\Beijing, 100048, P.R.China}
\end{center}

\begin{abstract}
 Given a map with underlying graph $\mathcal{G}$, if the set of  prime divisors of $|V(\mathcal{G}|$ is denoted by $\pi$, then we call the map a {\it $\pi$-map}.
 An  orientably-regular (resp. A regular )  $\pi$-map  is called {\it solvable}  if the group $G^+$ of all orientation-preserving automorphisms (resp. the group $G$ of  automorphisms)
   is solvable; and called {\it normal} if  $G^+$ (resp.  $G$) contains a normal $\pi$-Hall subgroup.

 In this paper, it will be proved that   orientably-regular $\pi$-maps are  solvable and normal if $2\notin \pi$  and  regular   $\pi$-maps are  solvable if $2\notin \pi$ and $G$ has no sections isomorphic to ${\rm PSL}(2,q)$ for some prime power $q$. In particular, it's shown that a regular   $\pi$-map with $2\notin \pi$ is normal if and only if $G/O_{2^{'}}(G)$ is isomorphic to a Sylow $2$-group of $G$.
 Moreover, nonnormal   $\pi$-maps  will be  characterized  and some properties   and constructions of   normal   $\pi$-maps  will be given in respective sections.
\end{abstract}


\renewcommand{\thefootnote}{\empty}
\footnotetext{$^*$corresponding author}
\footnotetext{ Email:  tianyao202108@163.com }

\section{Introduction}
A {\em map} is a cellular decomposition of a closed surface.
An alternative way to describe  maps is to consider them as cellular embeddings of graphs into closed surfaces.
By an {\em automorphism of a map} ${\cal M}$ we mean an  automorphism of the {\em underlying graph} ${\cal G}$
which extends to a   self-homeomorphism of the surface.
These automorphisms  form a subgroup $\Aut({\cal M})$ of the automorphism  group $\Aut({\cal G})$ of $\cal G$.
It is well-known that  $\Aut({\cal M})$ acts semi-regularly on the set of all flags (in most cases, which are incident vertex-edge-face triples).
If the action is regular, then we call the map as well as the corresponding embedding  {\em  regular}.

In the case of orientable supporting surface, if the group $\Aut^+(\MM)$ of all orientation-preserving
automorphisms of $\MM$ acts regularly on the set of darts (incident vertex-edge pairs) of $\MM$, then  we call $\MM$ an \emph{orientably-regular} map.
Such maps fall into two classes: those that admit also orientation-reversing automorphisms,
which are called \emph{reflexible}, and those that do not, which are  called \emph{chiral}. Therefore, a reflexible map is a regular map but a chiral map is not.

One of the central problems in topological graph theory is to classify all the regular  maps of  given underlying graphs.
In a general setting, the classification problem was treated in \cite{GNSS}.
However, the problem has been solved only for a few particular classes of graphs.
For instance,  complete graphs,  complete bipartite graphs, complete multipartite graphs,  $n$-dimensional cubes $Q_n$ and so on.
Moreover, there are closed connections between  regular maps and  other mathematical theories such as group theory, hyperbolic geometry and complex curves and so on.
During the past  forty  years, there have been abundant references on regular maps, see  \cite{JaJo,Jon2,Wil} for an overview.
\begin{defi}
A map is called a {\it $\pi$-map} if the set of prime divisors of number of vertices is denoted by $\pi$.
\end{defi}

  In this paper, we shall concentrate on   orientably-regular $\pi$-maps and  regular   $\pi$-maps.

A few attention has been paid to the situation that $\pi$ consists of a single prime. For instance, every orientably-regular map of a complete graph must be  a $p$-map, see \cite{JaJo}.
 Every $n$-dimensional hypercube $Q_{n}$ contains $2^n$ vertices,  whose  orientably-regular maps and nonorientable regular maps  were classified in \cite{CCDKNW} and \cite{KN1}, respectively. Every regular map of complete bipartite graph $K_{2^e,2^e}$ has $2^{e+1}$ vertices and a classification of  orientably-regular embeddings of this graph played an important role in a final classification  of $K_{n,n}$,  see \cite{DJKNS1,DJKNS2,Jon1,KK3}. But there are many maps that are not $p$-maps, for instance, the map of $K_{n,n}$ for $n$ not a $2$-power and the map of Petersen graph, etc. On the other hand, groups with suitable relations allow us to construct corresponding orientably-regular $\pi$-maps  and  regular $\pi$-maps, which means that the case $\pi$ contains more than one prime is ubiquitous from a group-theoretical point of view. For example, the simple group ${\rm PSL}(2,p)$ can be generated by two elements one of which is an involution or it can be generated by three involutions two of which commute, and we can use these relations to construct required orientably-regular $\pi$-maps  and  regular $\pi$-maps such that $\pi$ is not a single prime. In fact, if we require $\pi$ to contain only odd primes, then we will see that the group ${\rm PSL}(2,p)$ is involved in most of our situations when dealing with regular maps.

 One will see that the studies of {\it normal $\pi$-maps} is very closed to that of finite $\pi$-groups and its automorphisms as well. In other words, the studies of normal $\pi$-maps  stimulate us to focus on  some related $\pi$-groups and
 it might  pose some new research problems  for group theorists.

 During the studying of $p$-maps in \cite{DTL},  we obtained some new observations, ideas, methods and conjectures.
 This motivates us to pay more attention on  the more general $\pi$-maps. The aim of this paper is to give  a basic theoretical  characterization  for  orientably-regular $\pi$-maps and regular   $\pi$-maps.

  As usual, by  $\B{D}_k$,  $O_{\pi}(G)$ and $H\lhd G$,  we denote the dihedral group of order $k$, the maximal normal $\pi$-subgroup of $G$, and   a normal subgroup $H$ of  $G$,  respectively.
 First  we introduce   the following   concepts.

\begin{defi}  An  orientably-regular (resp. A regular )  $\pi$-map $\MM $ is called {\it solvable}  if  $\Aut^+(\MM)$  (resp. $\Aut(\MM)$) is solvable; and called {\it normal} if
 $\Aut^+(\MM)$ (resp.  $\Aut(\MM))$ contains a normal  $\pi$-Hall subgroup.
\end{defi}

\begin{rem}
  Suppose that $\MM$ is  an orientably-regular $\pi$-map which is reflexible. Then it is regular. Since $|\Aut(\MM):\Aut^+(\MM)|=2$, it follows that   $\Aut(\MM)$ is solvable if and only if
  $\Aut^+(\MM)$ is solvable. Take  a $\pi$-Hall subgroup $P$ of $\Aut(\MM)$.  If $2\notin \pi$,  then  $P\lhd  \Aut(\MM)$ if and only if $P\lhd \Aut^+(\MM)$.
   If $2\in \pi$,  then $P\lhd \Aut(\MM)$ implies that $P\cap \Aut^+(\MM)$ is the normal $\pi$-Hall subgroup of   $\Aut^+(\MM)$. However,  the converse is not true.
       For   example, consider the regular tetrahedron: the $2$-Hall subgroup of $\Aut^+(\MM)\cong A_4$ is normal   but every $2$-Hall subgroup of  $\Aut(\MM)\cong S_4$ is nonormal.
      \end{rem}

Now we are ready to state the main theorems of this paper.
\begin{theorem}\label{main1}  Let $\MM$ be an orientably-regular $\pi$-map  and $2\notin \pi$. Then
 \begin{enumerate}
   \item[{\rm (1)}] $\MM$ is solvable;
     \item[{\rm (2)}] $\MM$ is normal.
     \end{enumerate}
\end{theorem}

\begin{theorem}\label{main2}  Let $\MM$ be a regular $\pi$-map and $2\notin \pi$. Then
 \begin{enumerate}
   \item[{\rm (1)}] $\MM$ is solvable if $G$ has no sections isomorphic to ${\rm PSL}(2,q)$ for some prime power $q$;
     \item[{\rm (2)}] $\MM$ is normal if and only if $G/O_{2^{'}}(G)$ is isomorphic to a Sylow $2$-group of $G$.
     \end{enumerate}
\end{theorem}

\begin{rem}
  For a regular $\pi$-map $\MM$, it is known that $G={\rm Aut}(\MM)$ is generated by its involutions. If $2\in \pi$, then a normal $\pi$-Hall subgroup must be the whole group $G$ and hence $G$ itself is a $\pi$-group, which is of less interest to us. That is one of the reasons why we assume $2\notin \pi$ in Theorem 1.5.
      \end{rem}

It follows from Theorem~\ref{main1} and Theorem~\ref{main2} that   an orientably-regular  $\pi$-map can be  nonnormal  only if $2\in \pi$ and a  regular  $\pi$-map can be nonnormal only if $G/O_{2^{'}}(G)$ is not isomorphic to a Sylow $2$-group of $G$. To give a description in details for these  nonnormal $\pi$-maps, we
introduce three  families of  (multi)graphs and  maps as well.

   (i) By ${\cal D}_d$  we denote the $d$-dipoles,  which   is the graph of two vertices  joined by $d$ parallel edges. Every orientably-regular map of ${\cal D}_d$ is isomorphic to ${\cal D}(d, e):=\MM(G; x, x^ey)$ (see \cite{NS}), where $G=\lg x,y\di x^d=y^2=1, x^y=x^e\rg$ and $e^2\equiv 1(\mod m).$

   (ii)  By ${\cal S}_d$  we denote the $d$-semistars,   which   is the single vertex with $d$ semi-edges. By $DM(d)$ and $EM(d)$, we denote the regular map of  ${\cal S}_d$ in a disc and a sphere, respectively (see \cite{LS}). Note that in this paper  we only consider surfaces without boundary, except for the case in here we include the disc as a surface).

   (iii) By   $C_n^{(m)}$  we denote the graph resulting from the cycle $C_n$ of length $n$ by replacing each
edge with $m$ parallel edges.  There exists a nonorientable regular map  of $C_n^{(m)}$ if and only if $n=3,m=2$ and we denote the map by  ${\cal C}(3,2)$ (see \cite{HNSW}).

 Let $\mathcal{M}$ be a solvable orientably-regular map (resp. solvable regular map) and we denote by $\pi$ the set of prime divisors of the degree of the underlying graph of the map $\mathcal{M}$. If $G={\rm Aut}^+(\mathcal{M})$(resp. $G={\rm Aut}(\mathcal{M})$) does not contains a normal $\pi$-Hall subgroup, then we must have $2\in \pi$ (resp. either $2\in \pi$ or $3\in \pi$). In general, we have the following main theorem.

\begin{theorem}\label{main}  Suppose that $\MM$ is  a nonnormal  orientably-regular $\pi$-map  or a  nonnormal  regular   $\pi$-map. Let $G=\Aut^+(\MM)$ or $\Aut(\MM)$ and suppose that $G$ is solvable. Denote ${\overline \MM}$ by the quotient   map induced by $O_{\pi}(G)$, then one of the followings holds:
 \begin{enumerate}
     \item[{\rm (1)}] $2\in\pi$,   ${\overline \MM}={\cal D}(m,e)$, where $m\ge 3$ is odd and $e^2\equiv 1(\mod m)$ but $e\not\equiv 1(\mod m)$;
   \item[{\rm (2)}] $2\in \pi$,   ${\overline \MM}=DM(m)$ and  $\MM$ is nonorientable and nonnormal regular;
    \item[{\rm (3)}] $2\in \pi$,   ${\overline \MM}=EM(m)$ and  $\MM$ is normal  orientably-regular  but nonnormal regular;
     \item[{\rm (4)}] $2\notin \pi$ and $3\in \pi$,    ${\overline \MM}={\cal C}(3,2)$ and $\MM$ is   nonorientable and nonnormal  regular.
     \end{enumerate}
\end{theorem}

Theorems~\ref{main} will be proved by combining several lemmas in Sections 5.
\begin{rem}
  Note that it may happens that $G$ itself is a $\pi$-group. Then $G=P$, which is of less interest to us. So we concentrate on the case that $G$ is not a $\pi$-group, for instance $2\notin \pi$. For a normal orientably-regular $\pi$-map,   $G=P\!\rtimes \B{Z}_{m}$; and   for a normal regular $\pi$-map, $G=P\!\rtimes \B{D}_{2m}$. Furthermore, for such normal $\pi$-maps,   the studies  of  $\pi$-maps are essentially determinations of  automorphisms of a given $\pi$-group. Of course, it is  very difficult, because of  complexities of finite $\pi$-groups even for $\pi$ consisting of a single prime. So in general, we believe that it's hopeless to give a detailed description  of all normal $\pi$-maps. But in some sense, Theorem~\ref{main} is a  basic theoretical characterization of  orientably-regular $\pi$-maps and  regular $\pi$-maps, and it could  be a
starting point for studying such maps.
  \end{rem}

After this introductory section,  a brief introduction to regular maps and some   known results    will be given in Section 2; the solvability of  orientably-regular $\pi$-maps
and regular $\pi$-maps under suitable conditions will be proved in Section 3, some examples of $\pi$-maps are presented in the last of this section to explain why we focus on the case $2\notin \pi$; with some hypothesis, the normality of such maps will be discussed in Section 4; using the concept of quotient maps, the characterization of nonnormal solvable $\pi$-maps will be given in Section 5; the characterization of primitive $\pi$-maps is given in Section 6.

Except for the notations mentioned in the above,  more notations and terminologies     used in this paper are listed below.

\vskip 3mm

$|G|$, $|g|$:   the order of a group $G$ and the order of an element  $g$ in $G$, respectively;

 $|G:H|$,   $H\lhd G$ and $H\char G$:  the index of $H$ in $G$,    a normal subgroup  $H$ of $G$ and  a  characteristic subgroup $H$ of $G$, respectively;

$H_G$: the maximal normal subgroup of $G$ contained in the subgroup $H$;

$K\rtimes H$:  the semi-product of $K$ by $H$ where $K$ is normal;

$\pi '$: the complementary set of a subset   $\pi $ of   $\mathbb{P}$,  where $\mathbb{P}$ is the set of primes;

$\pi$-subgroup $H$:  a subgroup $H$ such that  every prime divisor of $|H|$  is contained in a subset $\pi $ of  $\mathbb{P}$;

$\pi$-Hall subgroup $H$: a $\pi$-subgroup $H$ such that   $(|H|, |G|/|H|)=1$;

$A-B$: those elements of $A$ not contained in $B$;

$\Syl_p(G)$:  the set of Sylow $p$-subgroups of $G$;

$F(G)$: the Fitting subgroup of $G$ (the product of all  nilpotent normal  subgroups of $G$);

 $O_{\pi}(G)$: the maximal normal $\pi$-subgroup of $G$ (the intersection of maximal $\pi$-subgroups of $G$).

\section{Preliminary Results}
In this section,  we shall   give a brief description  for  orientably-regualr and  regular maps and list some known  group theoretical   results used in this paper.

 \subsection{Regular maps  and orientably-regular maps}
(1) {\bf Regular Maps.}
\vskip 3mm
 A regular map can be described in the following way.
\begin{defi}\label{map1}
{\rm For a given finite set $F$ and three fixed-point-free involutory permutations $t, r, \ell $ on $F$, a quadruple
$\MM=\MM(F; t, r, \ell )$ is called a {\it combinatorial map} if they satisfy two conditions:
(1)\ $t\ell =\ell t$; (2)\ the group $\lg t,r, \ell \rg $ acts transitively on $F.$}
\end{defi}

For a given combinatorial map $\MM=\MM(F; t, r, \ell ),$  $F$  is called the {\it flag} set,
$t,  r, \ell$ are called {\it transversal, rotary, and longitudinal involution,} respectively.
The group $\lg t, r, \ell \rg $ is called the {\it monodromy group} of $\MM$, denoted by $\Mon(\MM)$.
We define the {\it vertices, edges} and {\it face-boundaries} of $\MM$ to be
respective the orbits of the subgroups $\lg t, r\rg$, $\lg t, \ell \rg $ and $\lg r, \ell \rg $ under the action on $F$,
incidence being given by a nontrivial set intersection.

The map $\MM$ in Definition~\ref{map1} is {\it unoriented}. Clearly, the even-word subgroup $\lg tr, r\ell \rg $
of $\Mon(\MM)$ has  index at most 2.  If the index is 1, then $\MM$ is said to be {\it nonorientable}.
If the index is 2, then one may fix an orientation for $\MM $ and so $\MM $ is said to be
{\it orientable}, while  the group $\Aut^+(\MM)$ of all orientation-preserving automorphisms
  acting transitively on all arcs is exactly $\lg tr, r\ell\rg .$

An automorphism of $\MM$ is defined to be a permutation on $F$ that commutes with $t,r$ and $\ell$, the automorphisms of $\MM $ form a group $\Aut (\MM)$  which is called the {\em automorphism group} of the map $\MM$. By definition, $\Aut(\MM )=C_{S_F}(\Mon (\MM)),$  the centralizer of $\Mon(\MM)$ in
$S_F$. It follows from the transitivity of $\Mon (\MM )$ on $F$ that   $\Aut (\MM )$ acts semi-regularly on $F.$  If the action is
regular, we call the map $\MM $ {\it regular}. As a consequence of a result in permutation group theory
(see \cite[I.Theorem 6.5]{Hup}), we get that for a regular map $\MM ,$ the two associated permutation groups $\Aut (\MM )$ and $\Mon (\MM )$
can be viewed as the right   regular  representations $R(G)$ and left regular representations  $L(G)$ of the abstract group $G\cong\Aut(\MM )\cong\Mon (\MM )$ mutually centralizing each other in $S_F$.
If there is no confusion we also identify  the elements $R(g)$ and $L(g)$  with  the elements $g\in G$,
 so that $\MM\cong \MM(G; t, r, \ell )$, which  is called an {\it algebraic map}.
Thus the subgroups $\lg t, r\rg$, $\lg t, \ell \rg $ and $\lg r,\ell \rg $, respectively, stand for the
stabilizer of a vertex, an edge and a face which are mutually incident.

\vskip 3mm (2) {\bf Orientably-regular maps.}
\vskip 3mm
An  orientably-regular map can be described in the following direct way.
 \begin{defi}\label{map2}
{\rm For a given finite set $D$ and two fixed-point-free  permutations $r, \ell $ on $D$ where $\ell $ is an involutory,  a triple
$\MM=\MM(D; r, \ell )$ is called a {\it combinatorial orientable   map} if $\lg r, \ell \rg $ acts transitively  on $D$.}
\end{defi}
For a given map $\MM=\MM(D; r, \ell ),$  $D$  is called the {\it darts},
$r, \ell$ are called {\it local rotation and arc-revision involution,} respectively.
The group $\lg r, \ell \rg $ is called the {\it monodromy group} of $\MM$, denoted by $\Mon(\MM)$.
We define the {\it vertices, edges} and {\it face-boundaries} of $\MM$ to be respective the
orbits of the cyclic subgroups $\lg r\rg$, $\lg \ell \rg $ and $\lg r\ell \rg $ under the action on $D$,
incidence being defined by a nontrivial set  intersection.

Similarly, one may define automorphisms of $\MM$ and know that
  $\Aut (\MM )$ acts semi-regularly on $D$.  If the action is
regular, we call the map $\MM $ {\it regular}. In this case,
  the two associated permutation groups $\Aut (\MM )$ and $\Mon (\MM )$
can be viewed as the right   regular  representations $R(G)$ and left regular representations  $L(G)$ of the abstract group $G=\lg r, \ell\rg \cong\Aut(\MM )\cong\Mon (\MM )$ mutually centralizing each other in $S_D$ so that the map is denoted by $\MM(G; r, \ell )$.

\subsection{Some Known Results }
\begin{prop}
\label{isa}
{\rm (\cite[Chap.5, Corollary 5.14]{ISA})}
Let $P\in Syl_p(G)$, where $G$ is a finite group and $p$ is the smallest prime divisor of $|G|$, and assume that $P$ is cyclic.
Then $G$ has a normal $p$-complement.
\end{prop}

\begin{prop}
\label{nc}
{\rm (\cite[Chap.1, Theorem 6.11]{SUZ})}
Let $H$ be a subgroup of a group $G$.
Then $C_G(H)$ is a normal subgroup of
$N_G(H)$ and the quotient $N_G(H)/C_G(H)$ is isomorphic
to a subgroup of $\Aut (H)$.
\end{prop}

\begin{prop}
\label{od}
{\rm (\cite{WJ})}
Every finite group of odd order is solvable.
\end{prop}

\begin{prop} \cite[Theorem 1]{GW} \label{dihedral} Let $G$ be a finite group with dihedral Sylow $2$-subgroups. Let
$O_{2'}(G)$ denote the maximal normal subgroup of odd order. Then $G/O_{2'}(G)$ is isomorphic to either a subgroup of ${\rm P\Gamma L}(2,q)$ containing $\PSL(2, q)$ where $q$ is odd,
or $A_7$, or a Sylow $2$-subgroup of $G$.
 \end{prop}

\begin{prop} \cite[Theorem 1.1]{MD} \label{generation} The alternating group $A_n$ can be generated by three involutions, two of which commute if and only if $n=5$ or $n\geq 9$.
 \end{prop}

\begin{prop} {\cite{Hup}} \label{sol}  Let $\{p_1, p_2, \cdots, p_l\}$ be the set of prime divisors of $|G|$. Suppose $G$ is solvable. Then
\begin{enumerate}\item[{\rm (i)}]
$F(G)=O_{p_1}(G)\times O_{p_2}(G)\times \cdots \times O_{p_l}(G)$;
\item[{\rm (ii)}] $C_G(F(G))\le F(G)$;
\item[{\rm (iii)}]  $\Phi(G)\le F(G)$.
\end{enumerate}
 \end{prop}

\section{Solvability  of $\pi$-Maps}
\begin{lem}
Let $G$ be a finite group with a cyclic subgroup of odd index $n$. Then $G$ is solvable.
\end{lem}
\demo By hypothesis, we know that $G$ has a cyclic Sylow-2 subgroup $S$. Thus by Proposition \ref{isa}, we know that $G$ has a normal subgroup $N$ of odd order such that $G=NS=N\rtimes S$. By Proposition \ref{od}, $N$ is solvable. Also $G/N\cong S$ is solvable, we derive that $G$ is solvable, as required.\qed

\begin{lem}
Suppose $G=\langle t, r, \ell\rangle$  such that $t,r, \ell$ are involutions, $t\ell =\ell t$ and $\langle r, t\rangle\cong D_{2k}$   where $|G:\langle r, t\rangle |$ is an odd number $n$. If $G$ has no sections isomorphic to ${\rm PSL}(2,q)$ for some prime power $q$, then $G$ is solvable.
\end{lem}
\demo  Under our condition, we know that either the order of a Sylow-2 subgroup of $G$ is no more than 2 or $G$ has a dihedral Sylow-2 subgroup. In the former case, we derive that $G=O_{2'}(G)S$ for a Sylow-2 subgroup of $G$, therefore $G$ is solvable by Proposition \ref{od}. Now we are left with the later case. By Proposition~\ref{dihedral},  we get $G/O_{2^{'}}(G)$ isomorphic to $A_{7}$, a group $N$ such that ${\rm PSL}(2,q)\leq N\leq  {\rm P\Gamma L}(2,q)$ with $q$ being an odd prime power or a Sylow $2$-subgroup of $G$. From Proposition ~\ref{generation}, we know that $A_7$ can not be generated by three involutions, two of which commute. Thus the first possibility is excluded. Now by assumption, we know that $G/O_{2^{'}}(G)$ must be isomorphic to a 2-group, hence the solvability of $G$ follows by Proposition \ref{od} and the fact that any $2$-group is solvable.\qed

\begin{lem} Suppose $G=\langle t, r, \ell\rangle$  such that $t,r, \ell$ are involutions, $t\ell =\ell t$ and $\langle r, t\rangle\cong D_{2k}$   where $|G:\langle r, t\rangle |$ is an odd number $n$ and we denote by $\pi$ the set of prime divisors of $n$. If $G$ is solvable, then $G/O_{2^{'}}(G)$ is isomorphic to  either $S_4$ or  a Sylow $2$-group of $G$ and $G/O_{2^{'}}(G)\not \cong S_4$ provided $3\notin \pi$. If $G$ is insolvable, then $G$ has a unique nonabelian composition series isomorphic to ${\rm PSL}(2,q)$ for some prime power $q$.
\end{lem}
\demo Suppose $4\nmid |G|$.   Then  $|G|=2m$ for some odd integer $m$ and so $G$ has a unique subgroup of order $m$, that is $O_{2'}(G)$. Clearly,  $G/O_{2'}(G)\cong \mathbb{Z}_2$,  as desired. Suppose $4| |G|$. Then $G$ contains a dihedral Sylow $2$-subgroup.  By Proposition~\ref{dihedral},  we get $G/O_{2^{'}}(G)$ isomorphic to $A_{7}$, a group $N$ such that ${\rm PSL}(2,q)\leq N\leq  {\rm P\Gamma L}(2,q)$ with $q$ being an odd prime power or a Sylow $2$-subgroup of $G$ while we exclude the possibility $A_7$ by Proposition ~\ref{generation} again.

 If $G$ is solvable. Since ${\rm PSL}(2,q)$ is simple when $q\ge 5$,  the first case was excluded and the second case  happens only if $N={\rm PSL}(2,3)\cong A_4$ or ${\rm P\Gamma L}(2,3)\cong S_4$. Note that $A_4$ cannot be generated by its involutions,  we get that  $G/O_{2^{'}}(G)$ is isomorphic to either  $S_4$ or a Sylow $2$-subgroup of $G$.
Now we assume that $3\notin \pi$.
Let $P$ be a $\pi$-Hall subgroup of $G$. Suppose that  $G/O_{2^{'}}(G)\cong S_4$.  From  $PO_{2'}(G)/O_{2'}(G)\lessapprox S_4$, we get that $P\leq O_{2'}(G)$ and so $G/O_{2'}(G)=\overline{P} \overline{H}=\overline{H} $, but $S_4$ can not be homomorphic image of a dihedral group, a contradiction.

If $G$ is insolvable, then $G/O_{2'}(G)$ is in the second case. Since ${\rm Aut}({\rm PSL}(2,q))={\rm PGL}(2,q)\rtimes {\rm Gal}(F_q/F_p)$, we have that $N/{\rm PSL}(2,q)$ is solvable, therefore the second statement follows from Jordan-H\"{o}lder Theorem and the fact that every composition series of a solvable group is of prime order.\qed

In the end of this section, we will give several examples to illustrate why we add the hypothesis $2\notin \pi$ in our study of orientably-regular $\pi$-maps and regular $\pi$-maps.

\begin{exam}
It is known that the symmetric group $S_7$ is generated by the two elements $r=(1~2\cdots~7)$ and $\ell=(1~2)$. Clearly we have that $\langle r\rangle$ is of even index in $G$. Calculation shows that the orientably-regular map $\MM(G,r,\ell)$ is reflexible.
\end{exam}

\begin{exam}
Let $G$ be the simple group ${\rm PSL}(2,8)$. Let $a$ be a generator of $F_8^*$ and $b$ an arbitrary nonzero element in $F_8$.
Set $r=\overline{\left(
                                                                 \begin{array}{cc}
                                                                   a & 0 \\
                                                                   -ab+a^{-1}b & a^{-1} \\
                                                                 \end{array}
                                                               \right)}
$ and $\ell=\overline{\left(
           \begin{array}{cc}
             1 & b \\
             0 & 1 \\
           \end{array}
         \right)}
$, then $r$ is of order $7$ and $\ell$ is an involution. Direct computation shows that $G=\langle r,\ell\rangle$. Set $g=\overline{\left(
                        \begin{array}{cc}
                          1 & b^{-1} \\
                          0 & 1 \\
                        \end{array}
                      \right)}$, then it is contained in  $N_G(\langle r\rangle)\cap C_G(\langle l\rangle)$ and therefore $r^g=r^{-1},\,\ell^g=\ell$. Thus the orientably-regular map $\MM(G,r,\ell)$ is reflexible.

\end{exam}
\begin{rem}
With the help of Magma, we find that the automorphism groups of chiral maps of low genus are all solvable. We does not know whether it holds for all chiral maps or not.
\end{rem}

\begin{exam}
Let $G=S_7$ be the symmetric group on seven letters, then $G$ is an insolvable group which has no sections isomorphic to any $2$-dimensional special linear group. Let $r=(2~7)(4~5)(3~6)$ and $t=(1~2)(3~7)(4~6)$ and $\ell=(1~2)$. Then we have that $G=\langle r,t,\ell\rangle$ and $\langle r,t\rangle$ is of even index in $G$, noting the regular map $\MM(G,r,t,\ell)$ is orientable.
\end{exam}

\begin{exam}
Let $G=S_8$ be the symmetric group on seven letters, then $G$ is an insolvable group which has no sections isomorphic to any $2$-dimensional special linear group. Let $r=(2~8)(3~7)(4~6)$ and $t=(1~2)(3~8)(4~7)(5~6)$ and $\ell=(1~2)$. Then we have that $G=\langle r,t,\ell\rangle$ and $\langle r,t\rangle$ is of even index in $G$, noting the regular map $\MM(G,r,t,\ell)$ is nonorientable.
\end{exam}

\section{Normality of $\pi$-Maps }

\begin{lem}
Let $\mathcal{M}$ be an orientably-regular map whose underlying graph has odd degree $n$, we denote by $\pi$ the set of prime divisors of $n$.  Then $G={\rm Aut}^+(\mathcal{M})$ has a normal $\pi$-Hall subgroup.
\end{lem}
\demo Since $\mathcal{M}$ is an orientably-regular map, $G=\langle r,\ell\rangle$ is solvable by Lemma 3.1. Let $H=\langle r\rangle$ and $P$ be a $\pi$-Hall subgroup of $G$ so that $G=PH$ and $|r|$ is even.  For the contrary,  assume that $P$ is not normal in $G$. Then $O_{\pi}(G)<P$.
Let $\overline{G}=G/O_{\pi}(G)$. Then  $O_{\pi}(\overline{G})=1$. Let $H_1$ be the ${\pi}'$-Hall subgroup of $H$.   Then $\overline{G} =\overline{P} \overline{H}=\overline{P}\overline{H_1}$.

   Since $O_{\pi}(\overline{G})=1$ and $|G|=n|H|$,  the Fitting subgroup  $F(\overline{G})$ of $\overline{G}$ is contained  in $\overline{H_1}$.  Since  $\overline{G}$ is solvable and $\overline{H_1}$  is cyclic,  it follows from Proposition \ref{sol} that $\overline{H_1}\le C_{\overline{G}}(F(\overline{G}))\leq F(\overline{G})$, which gives $F(\overline{G})=\overline{H_1}.$ In particular,
   $\overline{H_1}\lhd \overline{G}$, forcing  $C_{\overline{G}}(\overline{H_1})\unlhd \overline{G}$.

Since $\overline{H}\leq C_{\overline{G}}(\overline{H_1})$,  $C_{\overline{G}}(\overline{H_1})$ contains a Sylow $2$-subgroup of $\overline{G}$. Since $C_{\overline{G}}(\overline{H_1})\lhd \overline{G}$,   it contains all Sylow $2$-subgroups of $\overline{G}$ and
 then  contains all involutions of $\overline{G}$.  Now we have  $\overline{H}, \overline{\ell}\le C_{\overline{G}}(\overline{H_1}),$ which implies  $C_{\overline{G}}(\overline{H_1})=\overline{G}$.

     Finally, from $C_{\overline{G}}(\overline{H_1})=\overline{G}$, we get $[\overline{H_1}, \overline{P}]=\overline{1}$,  which implies   $\overline{P}\unlhd \overline{P}\overline{H_1}=\overline{G}$, that is $P\lhd G$, a contradiction.\qed

\f{\bf Proof of Theorem 1.4} Combining Lemma 3.1 and Lemma 4.1, we get the conclusion of Theorem 1.4.\qed

\vskip 3mm

\begin{exam} With notations as above, it can be shown that the automorphism group of a solvable orientably-regular map always has a normal $\pi_0$-Hall subgroup where $\pi_0=\pi-\{2\}$.
Note that Lemma $4.1$ is no longer true for $n$ even since we can not guarantee the solvability of $G$ in this situation. For instance, the group $G=\rm PSL(2,7)$ can be generated by two elements $r,\ell$ where $$r=\overline{\left(
                                                                                                                                          \begin{array}{cc}
                                                                                                                                            1 & 1 \\
                                                                                                                                            0 & 1 \\
                                                                                                                                          \end{array}
                                                                                                                                        \right)},~~~~\ell=\overline{\left(
                                                                                                                                                           \begin{array}{cc}
                                                                                                                                                             0 & -1 \\
                                                                                                                                                             1 & 0 \\
                                                                                                                                                           \end{array}
                                                                                                                                                         \right)}.
                                                                                                                                        $$
It is clear that the index of $\langle r\rangle$ in $G$ is even, but $G$ is simple which can not possess any nontrivial normal subgroups.\qed
\end{exam}

Now we turn to the determination of normality of a regular $\pi$-map. But it seems hopeless to give a complete answer for this question by our approach. As indicated by Lemma 3.3, $G/O_{2'}(G)\cong S_4$ or a Sylow-$2$ subgroup of $G$ if $G$ is solvable and $n$ is odd,  thus we will mainly deal with these two cases in the following.

\begin{lem} Let $\mathcal{M}$ be a regular map with $G={\rm Aut}(\mathcal{M})$ and we denote by $\pi$ the set of prime divisors of $|G:\langle r, t\rangle |$. If $G/O_{2^{'}}(G)$ is isomorphic to a Sylow $2$-group of $G$ , then $G$ has a normal $\pi_0$-Hall subgroup where $\pi_0=\pi-\{2\}$.
\end{lem}
\demo In particular, we know that $G$ is solvable by Proposition \ref{od}. Note that $G=\langle r,\ell,t \rangle$ and $H=\langle r,t\rangle$. By hypothesis, we have  $G=L\rtimes S$ where $L$ is a characteristic subgroup of odd order of $G$  and  $S\in Syl_2(G)$. Let $P$ be a $\pi_0$-Hall subgroup of $G$. Then $L=PH_1$, where $H_1$ is the $\{\pi,2\}'$-Hall subgroup of $H$ and clearly, $H_1$ is cyclic.   Since $G=L\rtimes S$ and $L\lhd G$, we conclude that $O_{\pi_0}(L)=O_{\pi_0}(G)$ noting that $G$ is solvable.

Assume that $H_1=1$.  Then  $P=L\lhd G$. So in what follows, we assume that $H_1\ne 1$.
 Let $\overline{G}=G/O_{\pi_0}(G)$. For the contrary, suppose that $P\ntriangleleft G$, that is $\overline{P}\ne \overline{1}$. Then we  first derive the following two facts:
\vskip 3mm
{\it Fact 1:  $F(\overline{L})={\overline H_1.}$ }
\vskip 3mm
Since $\overline{L}$ is normal in $\overline{G}$, we get  $O_{\pi_0}(\overline{L})=1$, which in turn implies the Fitting subgroup $F(\overline{L} )$ of $\overline{L}$  is contained in  $O_{\pi'}(\overline{L} )$, that is the intersection of all $\pi'$-Hall subgroups of $\overline{L} $.
 Since  $H_1$ is a $\pi'$-Hall subgroup of $L$ and also a $\{\pi,2\}'$-Hall subgroup of $H$,  we have $O_{\pi'}(\overline{L} )\leq \overline{H}_1$.
   Since $\overline{G} $ is solvable and $\overline{H}_1$  is cyclic, we have  from Proposition~\ref{sol}.(ii) that $$O_{\pi'}(\overline{L})\le \overline{H}_1\le C_{\overline{L} }(O_{\pi'}(\overline{L}))\le  C_{\overline{L}}(F(\overline{L}))\leq F(\overline{L})\le O_{\pi'}(\overline{L}),$$
   which forces   $\overline{H}_1=F(\overline{L})\lhd \overline{G}.$
\vskip 3mm
{\it Fact 2: There exists a $\overline{Q}\in Syl_q(F(\overline{L}))$ for some $q\notin \pi_0$  such that $ |\overline{G}:C_{\overline{G}}(\overline{Q})|$ } is divided by twice a $\pi_0$-number.
\vskip 3mm
Since  $C_{\overline{L}}(F(\overline{L}))\le F(\overline{L})=\overline{H}_1$ (Fact 1) and $\overline{P} \varsubsetneq \overline{H}_1$,  we get  $[\overline{P}, F(\overline{L})]\ne \overline{1}$.
 Therefore, there exists some $\overline{Q}\in Syl_{q}(F(\overline{L}))$  such that $[\overline{P}, \overline{Q}]\ne \overline{1}$ where $q\notin \pi_0$. Since $\overline{Q}$ char $F(\overline{L})$ char $\overline{L}\lhd \overline{G}$,
  we have $\overline{Q}\unlhd\overline{G}$ and thus $C_{\overline{G}}(\overline{Q})\unlhd\overline{G}$.  Therefore,   $C_{\overline{G}}(\overline{Q}) $ contains one of $\pi_0$-Hall subgroups (resp. Sylow 2-subgroups) if and only if  it contains all $\pi_0$-Hall subgroups  (resp. all involutions).
   Note that $\overline{H}\cong D_{2n}$ has an involution $\overline{h} \in \overline{H}$ such that  $[\overline{Q}, \overline{h}]\ne 1$. Now  $\overline{P} $ and $\langle \overline{h} \rangle $  are not contained in  $C_{\overline{G}}(\overline{Q})$. Therefore, $C_{\overline{G}}(\overline{Q})$   contains neither a $\pi_0$-Hall subgroup of $G$ nor a Sylow
   2-subgroup of $G$, as required.

\vskip 3mm
Come back to our proof. By Proposition~\ref{nc}, we have
$$N_{\overline{G}}(\overline{Q})/C_{\overline{G}}(\overline{Q})=\overline{G}/C_{\overline{G}}(\overline{Q})\lessapprox \rm Aut(\overline{Q}).$$
 Since $\rm Aut(\overline{Q})$ is cyclic, we get that $\overline{G}/C_{\overline{G}}(\overline{Q})$ is cyclic. Let $\overline{T} /C_{\overline{G}}(\overline{Q})$ be the Sylow 2-subgroup of  $\overline{G}/C_{\overline{G}}(\overline{Q})$. By Fact 2,  $\overline{T}\neq \overline{1}, \overline{G}$.
 Since $\overline{T} $ is a subgroup of $\overline{G}$ of odd index,   $\overline{T} $ contains a Sylow $2$-subgroup of $\overline{G} $. Since $\overline{T}\lhd \overline{G} $, it contains all Sylow $2$-subgroups of $\overline{G}$, and in particular it contains all involutions.
   Then $\overline{G} =\langle \overline{r}, \overline{t}, \overline{\ell}\rangle \le \overline{T}$,  a contradiction.

   In summary,  we get $P\lhd G$. \qed
\vskip 3mm

\f{\bf Proof of Theorem 1.5} The first assertion follows by Lemma 3.3. To prove the second statement, it suffices to prove the ``only if'' direction by Lemma 4.3. Since $2\notin \pi$, Lemma 3.3 applies, we know that $G/O_{2'}(G)$ is isomorphic to either a Sylow-$2$ subgroup of $G$ or an almost simple group with socle ${\rm PSL}(2,q)$. We claim that $G$ can not possess a normal $\pi$-Hall subgroup in the latter case. Now assume that $G/O_{2'}(G)$ is isomorphic to an almost simple group with socle ${\rm PSL}(2,q)$, in particular, it can be neither cyclic nor dihedral. Since $G=HP$, we have $H/H\cap P\cong G/P$. Since $2\notin \pi$, $P$ is a normal subgroup of odd order and then is contained in $O_{2'}(G)$ by definition. Thus we have a surjective homomorphism from $G/P$ to $G/O_{2'}(G)$ and therefore a surjective homomorphism from $H$ to $G/O_{2'}(G)$. Since $H$ is dihedral, $G/O_{2'}(G)$ must be either cyclic or dihedral, which is a contradiction to our assumption, the proof is now completed.\qed

\vskip 3mm

\begin{exam}
With the same notations and hypothesis as in Lemma 4.3, there are situations showing that $G$ may not have a normal $\pi_0$-Hall subgroup if we only suppose that $G$ is solvable. The simplest example is $G=S_4=\langle r,t,\ell\rangle$ where $r=(13), \ t=(12)(34)$ and $\ell =(12)$, in which case $\pi=\pi_0=\{3\}$. But it's known that $G$ has no normal Sylow-$3$ subgroup.
\end{exam}

\begin{lem} With the same notations and hypothesis as in Lemma 4.3, if $G={\rm Aut}(\mathcal{M})$ satisfies that $G/O_{2'}(G)\cong S_4$, then $\MM$ is nonnormal and $G$ is a $\{\pi,2\}$-group.
\end{lem}
\demo We first claim that $2\notin \pi$ and $3\in \pi$ under our conditions. To see it, we write $\tilde{G}=G/O_{2'}(G)$ and $L$ for $O_{2'}(G)$. Then $\tilde{H}$ is a dihedral subgroup and $\tilde{G}=\langle \tilde{H}, \tilde{\ell}\rangle=\langle \tilde{r},\tilde{t},\tilde{\ell}\rangle.$ Since $S_4$ has trivial centre and it is not isomorphic to a dihedral group, we conclude that $\tilde{H}\cong D_8$, thus $3$ divides $|\tilde{G}:\tilde{H}|$. Since $|\tilde{G}:\tilde{H}|$ divides $|G:H|$, we conclude that $3\in \pi$. Since $O_{2'}(G)$ is a subgroup of odd order, that $2$ does not divide $|\tilde{G}:\tilde{H}|$ implies that $|G:H|$ is odd, as required. If $G$ has a normal $\pi$-Hall subgroup, then it must be contained in $O_{2'}(G)$ which contradicts that $3$ divides the order of $G/O_{2'}(G)$. Hence $G$ does not contain a normal $\pi$-Hall subgroup.

Next we are going to see that all prime divisors of the order of $G$ lie in $\{\pi,2\}$. Let $S$ be a Sylow-2 subgroup of $H$ and $P$ a $\pi$-Hall subgroup of $G$. Let $H_1$ be the $\{\pi,2\}'$-Hall subgroup of $H$, clearly $H_1$ is cyclic, we claim that $H_1=1$.

Suppose  $H_1\ne 1$. Then there exists a prime divisor $q$  of $|G|$  where $q\ne 2, 3$. Set $\overline{G}=G/O_{\pi}(G)$.
 Then we first show the following two facts.

\vskip 3mm
{\it Fact 1:  $F(\overline{L})=\overline {H_1}.$ }
\vskip 3mm
  Since $\overline{L}$ is normal in $\overline{G}$, we get  $O_{\pi}(\overline{L})=1$, which in turn implies the Fitting subgroup $F(\overline{L} )$ of $\overline{L}$  is contained in  $O_{\pi'}(\overline{L} )$, that is the intersection of all $\pi'$-Hall subgroups of $\overline{L} $.
 Since  $H_1$ is also a $\{\pi,2\}'$-Hall subgroup of $G$ and the index of $L$ in $G$ is $2^3\cdot 3$, we have $H_1\leq L$. Thus $H_1$ is a $\pi'$-Hall subgroup of $L$ and then $\overline{H_1}$ is a $\pi'$-Hall subgroup of $\overline{L}$, which implies that $O_{\pi'}(\overline{L} )\leq \overline{H}_1$.
   Since $\overline{G} $ is solvable and $\overline{H}_1$  is cyclic, we have  from Proposition~\ref{sol}.(ii) that $$O_{\pi'}(\overline{L})\le \overline{H}_1\le C_{\overline{L} }(O_{\pi'}(\overline{L}))\le  C_{\overline{L}}(F(\overline{L}))\leq F(\overline{L})\le O_{\pi'}(\overline{L}),$$
   which forces   $\overline{H}_1=F(\overline{L})\lhd \overline{G}.$
\vskip 3mm
 \vskip 3mm
{\it Fact 2: For any  $\overline{Q}\in Syl_q(F(\overline{L}))$ where $q\notin \{2, \pi\}$, we have that $6\mid |\overline{G}:C_{\overline{G}}(\overline{Q})|$. }
\vskip 3mm
 By Fact 1, $F(\overline{L})=\overline{H_1}$.   Let $\overline{Q}\in Syl_q(F(\overline{L}))$, where we may take $Q\le H_1\le H$.  Since $\overline{Q}~~{\rm char}\, F(\overline{L})~~{\rm char} \,\overline{L}$,  we have $\overline{Q}\unlhd \overline{G}$, forcing  $\overline{C}:=C_{\overline{G}}(\overline{Q})\unlhd \overline{G}.$ Since $[\overline{Q},\overline{S}]\neq \overline{1}$, some element of $\overline{S}\cong D_8$ is not contained in $\overline{C}$.  By the normality of $\overline{C}$, we know that $\overline{C}$ cannot contain any Sylow 2-subgroup of $\overline{G}$, that is $2\mid |\overline{G}:\overline{C}|$. But from $\overline{Q}\leq \overline{H}$, we know that a Sylow $2$-subgroup $\overline{A}$ of $\overline{C}$ is isomorphic to $\mathbb{Z}_4$.

 Suppose that $\overline{C}$ contains a Sylow 3-subgroup of $\overline{G}$. Then $3\mid |\overline{C}\overline{L}/\overline{L}|$.
  Since  $\overline{C}\overline{L}/\overline{L} \lhd \overline{G}/\overline{L}\cong G/L\cong S_4$, it follows that  $\overline{C}\overline{L}/\overline{L}$ contains a subgroup isomorphic to $A_4$,   contradicting to $\overline{A}\cong \mathbb{Z}_4$.
      Therefore,  $\overline{C}$  does not contain any Sylow 3-subgroup of $\overline{G}$,  that is $3\mid |\overline{G}:\overline{C}|$.

\vskip 3mm
Finally, using the same arguments as  that in  last paragraph of  Lemma 4.3, we may get  $P\lhd G$, a contradiction. Therefore, $H_1$ must be trivial, as required.\qed

\vskip 3mm
\section{Nonnormal $\pi$-Maps }
By the results in last section, it is now meaningful to characterize the solvable regular map(resp. solvable orientably-regular map) $\mathcal{M}$ whose automorphism group ${\rm Aut}(\mathcal{M})$(resp. ${\rm Aut}^+(\mathcal{M})$) does not possess a normal $\pi$-Hall subgroup, where $\pi$ is the set of prime divisors of the degree of the underlying graph of $\mathcal{M}$.

\begin{lem} Let $\mathcal{M}$ be a solvable orientably-regular map and we denote the set of prime divisors of the degree of the underlying graph of the map $\mathcal{M}$ by $\pi$ and $\overline {\mathcal{M}}$ the quotient map of $\mathcal{M}$ induced by $O_{\pi}(G)$. If $G={\rm Aut}^+(\mathcal{M})$ does not contain a normal $\pi$-Hall subgroup, then $\overline {\mathcal{M}}={\mathcal D}(m,e)$ and $m\ge 3$ is odd and $e^2\equiv 1(\mod m)$ but $e\not\equiv 1(\mod m)$.
\end{lem}
\demo Note that 2 must be a $\pi$-number by the results in Section 4. Set  $G=\langle r, \ell\rangle $ and  $\overline{G} =G/O_{\pi}(G)$. Let $P$ be a $\pi$-Hall subgroup of $G$. By hypothesis, we have $O_{\pi}(G)<P$. Then $\overline{G} =\overline{P} \langle \overline{r}\rangle $, and the Fitting subgroup $F(\overline{G})$ is of odd order, which implies  $F(\overline{G})\le  \langle \overline{r}\rangle$.
Since $\overline{G} $ is solvable,  we have  $\langle \overline{r}\rangle \le C_{\overline{G}}(F(\overline{G}))\le F(\overline{G})$, that is  $F(\overline{G})=\langle \overline{r}\rangle$, and in particular  $|\overline{r}|$ is odd.  Therefore,
$\overline{G}=\langle \overline{r},\overline{\ell}\rangle =\langle \overline{r}\rangle\langle\overline{\ell}\rangle$. Since $\overline{P}\ne \overline{1}$,  we get $\overline{\ell}\ne \overline{1}$  and so  $|P:O_{\pi}(G)|=2$. This   implies  that the quotient map ${\overline {\mathcal{M}}}$ of $\mathcal{M}$ induced by $O_{\pi}(G)$  has two vertices
and then the underlying graph of ${\overline {\mathcal{M}}}$  is  a dipoles ${\cal D}_m$, where $m=|\overline{r}|$ is odd.
A classification of orientably-regular maps of ${\mathcal{D}}_m$ was  given in \cite{NS}. As mentioned in Section 1,  every such map is isomorphic to  ${\mathcal{D}}(m,e)$, where $e^2\equiv 1(\mod m)$.
  Since $\overline{G}$ is nonabelian, $e\not\equiv 1(\mod m)$.   Therefore, $\mathcal{M}$ is a regular cover of ${\mathcal{D}}(m,e)$,   with the covering transformation group $O_{\pi}(G)$, a $\pi$-group.\qed

\vskip 3mm

\begin{lem} Let $\mathcal{M}$ be a solvable regular map and we denote by $\pi$ the set of prime divisors of the degree of the underlying graph of the map $\mathcal{M}$. Let ${\overline {\mathcal{M}}}$ be the quotient map of $\MM$ induced by the normal subgroup $O_{\pi}(G)$. If $G={\rm Aut}(\mathcal{M})$ does not contain a normal $\pi$-Hall subgroup, then either $2\in \pi$ or $3\in \pi$, and we have one of the following two situations:

(1)$2\in \pi$, in this case, ${\overline {\mathcal{M}}}=DM(m)$ and $\mathcal{M}$ is nonorientable and nonnormal regular or ${\overline {\mathcal{M}}}=EM(m)$ and $\mathcal{M}$ is normal orientably-regular but nonnormal regular;

(2)$2\notin \pi$ and $3\in \pi$, in this case, $\overline{\mathcal{M}}\cong {\cal C}(3,2)$.
\end{lem}
\demo First we deduce that either $2\in \pi$ or $3\in \pi$ by the results in Section 4.
\vskip 3mm
{\it Situation 1: $2\in \pi$.}
\vskip 3mm
Now if we are in the first situation, set $G=\langle r, t, \ell\rangle $   and $H_1$  the $\pi'$-Hall subgroup of $\langle r, t\rangle $.
 Set $\overline{G} =G/O_{\pi}(G)$ and let $P$ be a $\pi$-Hall subgroup of $G$. Then $\overline{G} =\overline{P}\langle \overline{r}, \overline{t}\rangle $. Then $|F(\overline{G})|$ is odd and so $F(\overline{G})\le \overline{H_1}$.   From $\langle \overline{r}\overline{t}\rangle \le C_{\overline{G}}(F(\overline{G}))\le F(\overline{G})$, that is  $F(\overline{G})=\langle \overline{r}\overline{t}\rangle \lhd \overline{G}$, and thus $m:=|\langle \overline{r}\overline{t}\rangle|$ is odd and  $\langle \overline{r}, \overline{t}\rangle \cong D_{2m}$. Therefore,
$\overline{G}=\langle \overline{r} \overline{t}\rangle \rtimes \langle \overline{t}, \overline{\ell}\rangle =\langle \overline{r}, \overline{t} \rangle \langle \overline{\ell}\rangle $, reminding $[t,\ell]=1$.

 Suppose that $\overline{\ell}\not\in \langle\overline{r}, \overline{t}\rangle $. Then
  $\overline {\mathcal{M}}$  has two vertices and so its underlying graph  is a dipoles ${\mathcal D}_m$ where
$m$ is odd.   If  $\mathcal{M}$ is nonorientable,  then   $\overline {\mathcal{M}}$ is nonorientable too. However, by  \cite[Lemma 3.1]{HNSW}, the  dipoles ${\mathcal D}_m$ has a nonorientable map if and only if $m=2$. Since $m$ is odd,
our  $\mathcal{M}$ is  orientable and (so reflexible).  Since $\overline{G} \cong \mathbb{Z}_m\rtimes D_4$ and $F(\overline{G})=C_{\overline{G}}(F(\overline{G}))\cong \mathbb{Z}_m$, clearly $m$ must contains at least two prime divisors. Clearly, $\overline{\langle rt, t\ell\rangle}$ has index 2 in $\overline{G}$ and so  $O_{\pi}(G)\le \langle rt, t\ell\rangle $.
  Therefore,  $\mathcal{M}$ is nonnormal orientably-regular if and only if it is nonnormal regular.

 Suppose that  $\overline{\ell}\in \langle\overline{r}, \overline{t}\rangle $. Then $\overline{G} \cong D_{2m}$ and so the quotient map $\overline {\mathcal{M}}$  has only one vertex. Since $\overline{\ell} \overline{t}=\overline{t}\overline{\ell}$ and $\overline{G} \cong D_{2m}$ where $m$ is odd, we have either $\overline{\ell}=\overline{1}$ or $\overline{\ell}=\overline{t}$:
  \vskip 3mm
 (1) $\overline{\ell}=1$: Following  \cite{LS}, $\overline {\mathcal{M}}$ is said to be degenerate. The underlying graph is a $m$-semistar and it can be embedded in a disc, which is an orientable surface with  boundary. This map is denoted by $DM(m)$. Since $\langle \overline{r}\overline{t}, \overline{t}\overline{\ell}\rangle =\langle \overline{r}, \overline{t}\rangle =\overline{G}$, we  have $\langle rt, t\ell\rangle =G$ and so $\mathcal{M}$ is nonorientable.

\vskip 3mm
 (2) $\overline{\ell}=\overline{t}$:  Following  \cite{LS}, $\overline {\mathcal{M}}$ is said to be redundant.  The underlying graph is a $m$-semistar too  and it can be embedded in a sphere.  This map is denoted by $EM(m)$. Since $\langle \overline{r}\overline{t}, \overline{t}\overline{\ell}\rangle =\langle \overline{r}\overline{t}\rangle \ne \overline{G}$, we  have $\langle rt, t\ell\rangle \ne G$ and $O_{\pi}(G)\le \langle rt, t\ell\rangle $  and so $\mathcal{M}$ is orientable and so reflexible. Moreover, $O_{\pi}(G)\le \langle rt, t\ell\rangle$, which  means
 that $\mathcal{M}$ is  normal orientably-regular  but nonnormal regular.
\vskip 3mm
{\it Situation 2: $2\notin \pi$ and $3\in \pi$.}
\vskip 3mm
 In this situation, if $G/O_{2'}(G)$ is isomorphic to a Sylow 2-subgroup of $G$,  then $\mathcal{M}$ is normal by Lemma 4.3.
Following Lemma 3.3, we need to consider the case that  $G/O_{2'}(G)\cong S_4$.

By our hypothesis, we know that $\overline{G}=G/O_{\pi}(G)\cong S_4.$ Note that $\langle\overline{r}\overline{t},\overline{t}\overline{\ell}\rangle$ is a normal subgroup of $\overline{G}$ containing a subgroup $\langle\overline{r}\overline{t}\rangle\cong \mathbb{Z}_4,$ we must have $\overline{G}=\langle\overline{r}\overline{t},\overline{t}\overline{\ell}\rangle$. On the other hand, that $|G:\langle rt,t\ell\rangle|\leq 2$ yields $O_{\pi}(G)\leq \langle rt,t\ell\rangle$, forcing $G=\langle rt,t\ell\rangle$. Hence $\mathcal{M}$ is a nonorientable map.

   Since $\overline{G} =\langle \overline{r}, \overline{t}, \overline{\ell}\rangle \cong S_4$ where $\overline{H}=\langle \overline{r}, \overline{t}\rangle \cong D_8$, $\overline {\mathcal{M}}$ has three vertices. Since the kernel of the action of $\overline{G}$ on vertices  is isomorphic to $D_4$, there are double edges (four flags) between every two adjacent vertices  and thus the underlying graph is $C_3^{(2)}$ and $\overline{\mathcal{M}}\cong {\cal C}(3,2)$. Therefore, $\mathcal{M}$ is a regular cover of ${\cal C}(3,2)$,   with the covering transformation group $O_{\pi}(G)$, a $\pi$-group.\qed

\f{\bf Proof of Theorem 1.7} Combining Lemma 5.1 and Lemma 5.2, we get the conclusions in Theorem 1.7.\qed

\section{primitive $\pi$-Maps }

\begin{theorem} Let $\mathcal{M}$ be an  orientably-regular map and we denote by $\pi$ the set of prime divisors of the degree of the underlying graph of the map $\mathcal{M}$. If $G={\rm Aut}^+(\mathcal{M}))$ acts primitively on the vertices, then $\pi$ consists of a single prime.
\end{theorem}
\demo First, we claim that if $G$ is a finite group with a cyclic maximal subgroup $H$, then $G$ is solvable. To see it, we use induction on the order of $G$. First we may assume that $H$ is core-free, otherwise $G/H_G$ would be solvable by inductive hypothesis, this together with the solvability of $H_G$ imply that $G$ is solvable and we will be done. In particular, $H$ has the trivial-intersection property, i.e., $H\cap H^g=1,~\forall g\in G-H$. Thus the well-known Frobenius Theorem \cite[Theorem 7.2]{IS} implies that $G$ is a Frobenius group with Frobenius complement $H$ and Frobenius kernel $N$. Since $G/N\cong H$ is solvable and $N$ is nilpotent by Thompson's thesis in \cite{JT}, we conclude that $G$ is solvable, as claimed. Since it's known that any maximal subgroup of a solvable group is of prime-power index \cite[Problem 3B.1]{ISA}, we deduce that $\pi$ must contain only one prime, as required.
\qed

\begin{theorem} Let $\mathcal{M}$ be a regular map and we denote by $\pi$ the set of prime divisors of the degree of the underlying graph of the map $\mathcal{M}$. If $G={\rm Aut}(\mathcal{M})$ acts primitively on the vertices and $2\notin \pi$ and we denote by $H$ a point stabilizer, then either $\pi$ consists of a single prime or $G/H_G$ is an almost simple group with socle ${\rm PSL}(2,q)$. In the latter case, $H/H_G$ is a $\pi'$-Hall subgroup of $G/H_G$ and the degree of the graph is either $q(q-1)/2$ or $q(q+1)/2$.
\end{theorem}
\demo If $2\notin \pi$, then $\pi$ consists of odd primes. Recall that $G=\langle r,t,\ell\rangle$ and $H=\langle r,t\rangle$ by our convention where $r,t,\ell$ are three involutions and $t\ell=\ell t$. If $G$ is solvable, then it is clearly that the index of $H$ in $G$ is a prime power, i.e., there is only one prime in $\pi$. Thus we may assume that $G$ is insolvable. Set $\overline{G}=G/H_G$. Then $\overline{H}$ is core-free and $\overline{G}$ is insolvable.  We know that $\overline{G}$ has a unique nonabelian composition factor isomorphic to ${\rm PSL}(2,q)$ for some prime power $q$ by Lemma 3.3. Now we can find a minimal normal subgroup $\overline{N}$ of $\overline{G}$ which is not contained in $\overline{H}$ by our assumption. Thus $\overline{N}$ has at most one composition series isomorphic to ${\rm PSL}(2,q)$. But on the other hand, it is known that $\overline{N}$ is the direct product of isomorphic simple groups. Hence either $\overline{N}\cong {\rm PSL}(2,q)$ or $\overline{N}$ is an elementary abelian $p$-group. In the latter case, $\overline{G}$ is solvable which contradicts our assumption. Now we are going to deal with the former case. As $C_{\overline{G}}(\overline{N})\lhd \overline{G}$ and $Z(\overline{N})=1$, we claim that $C_{\overline{G}}(\overline{N})=\overline{1}$ under our conditions. Suppose for the contrary $C_{\overline{G}}(\overline{N})\neq \overline{1}$, then it must contain some minimal normal subgroup $\overline{L}$ of $\overline{G}$. But the previous analysis implies that $\overline{L}\cong {\rm PSL}(2,q)$ again, which contradicts that $C_{\overline{G}}(\overline{N})\cong {C_{\overline{G}}(\overline{N}) \overline{N}/\overline{N}}\lhd \overline{G}/\overline{N}$ has no composition factors isomorphic to ${\rm PSL}(2,q)$, as required. Since $\overline{G}/C_{\overline{G}}(\overline{N})\hookrightarrow \rm Aut(\overline{N})$, we may regard $\overline{G}$ as an almost simple group, that is, ${\rm PSL}(2,q)\leq \overline{G}\leq {\rm P\Gamma L}(2,q)$. Now let $r$ be a prime in $\pi$ and $\overline{Q}$ be a Sylow-$r$ subgroup of $\overline{H}$. Then Sylow Theorem implies that $\overline{Q}$ is contained in some Sylow-$r$ subgroup $\overline{R}$ of $\overline{G}$. Since $\overline{Q}$ is normal in $\overline{H}$ and $N_{\overline{R}}(\overline{Q})>\overline{Q}$, we derive that $\overline{Q}\lhd \overline{G}$. But $\overline{G}$ is an almost simple group, which yields that $\overline{Q}=\overline{1}$, whence $\overline{H}$ is a $\pi'$-Hall subgroup.  Using the results on maximal subgroups of almost simple group with socle ${\rm PSL}(2,q)$ from \cite[Theorem 1.1]{MG}, we know that $\overline{H}\cap \overline{N}\cong D_{q-1}$ or $D_{q+1}$, therefore the degree of the graph is either $q(q-1)/2$ or $q(q+1)/2$, the proof is now completed.

\begin{rem}
If $2\in \pi$, much less could be said about $G$. For instance, many simple groups such as the Suzuki simple groups may contain a dihedral subgroup as a maximal subgroup and we can use this fact to construct primitive regular $\pi$-maps such that the conclusion in Theorem 6.2 is no longer true. The case $G$ is an almost simple group can really happen. It's known that except very few situations, $G={\rm PGL}(2,p)$ has a maximal subgroup isomorphic to $D_{2(p+1)}$ which is of index $p(p-1)/2$ in ${\rm PGL}(2,p)$. If we require $p\equiv 3 (mod~4)$, then $p(p-1)/2$ is an odd number. Suppose $D_{2(p+1)}=\langle r,t\rangle$ for two involutions $r,t$. Since $C_G(t)$ is also a dihedral group which is different from $\langle r,t\rangle$, we can choose some involution $\ell\in C_G(t)-\langle r,t\rangle$ so that $G=\langle r,t,\ell\rangle$ and $t,\ell$ commute.
\end{rem}

\f{\bf Note} The above two theorems show that in most cases primitive orientably-regular maps and primitive regular maps are $p$-maps. For primitive orientably-regular $p$-maps and primitive regular $p$-maps with $p\geq 5$, a characterization of their automorphism groups with emphasis on Sylow-$p$ subgroups can be found in \cite{DTL}.

\end{document}